\documentclass[12pt]{article}
\usepackage{amsmath,amsthm,amsfonts,amssymb,amscd}
\def\qed{{\unskip\nobreak\hfil\penalty50
\hskip2em\hbox{}\nobreak\hfil$\square$
\parfillskip=0pt \finalhyphendemerits=0\par}\medskip}
\def\proof{\trivlist \item[\hskip \labelsep{\bf Proof\ }]}
\def\endproof{\null\hfill\qed\endtrivlist}

\def\lan{\langle}
\def\ran{\rangle}
\def\Ad{{\mathrm {Ad}}}

\def\a{\alpha}

\def\e{\varepsilon}

\def\la{\lambda}

\def\phi{\varphi}

\def\Om{\Omega}

\def\r{{\rho}}
\def\t{{\tau}}

\def\emptyset{\varnothing}
\def\setminus{\smallsetminus}

\def\Diff{{\mathrm {Diff}}}
\def\Mob{{\rm\textsf{M\"ob}}}

\def\lan{\langle}
\def\ran{\rangle}
\def\Ad{{\mathrm {Ad}}}

\def\a{\alpha}

\def\e{\varepsilon}

\def\l{\lambda}

\def\phi{\varphi}

\def\Om{\Omega}

\def\r{{\rho}}
\def\t{{\tau}}

\newtheorem{theorem}{Theorem}[section]
\newtheorem{lemma}[theorem]{Lemma}
\newtheorem{conjecture}[theorem]{Conjecture}
\newtheorem{corollary}[theorem]{Corollary}
\newtheorem{definition}[theorem]{Definition}

\newtheorem{proposition}[theorem]{Proposition}
\newtheorem{remark}[theorem]{Remark}

\def\emptyset{\varnothing}
\def\setminus{\smallsetminus}

\def\Diff{{\mathrm {Diff}}}
\def\Mob{{\rm\textsf{M\"ob}}}

\def\res{\!\restriction\!}
\def\A{{\cal A}}

\def\B{{\cal B}}

\def\I{{\cal I}}

\def\L{{\cal L}}

\def\H{{\cal H}}

\def\Z{{\mathbb Z}}
\def\R{{\mathbb R}}

\def\1#1{{\bf #1}}
\def\2#1{{\mathcal #1}}
\def\3#1{{\sl #1}}
\def\4#1{{\tt #1}}
\def\5#1{{\sf #1}}
\def\6#1{{\mathfrak #1}}
\def\7#1{{\mathbb #1}}

\renewcommand{\qed}{\ \hfill $\blacksquare$}

\newcommand{\bdefin}{\begin{definition}}
\newcommand{\blemma}{\begin{lemma}}
\newcommand{\bprop}{\begin{proposition}}
\newcommand{\btheor}{\begin{theorem}}
\newcommand{\bcoro}{\begin{corollary}}
\newcommand{\bconj}{\begin{conjecture}}
\newcommand{\edefin}{\end{definition}}
\newcommand{\elemma}{\end{lemma}}
\newcommand{\eprop}{\end{proposition}}
\newcommand{\etheor}{\end{theorem}}
\newcommand{\ecoro}{\end{corollary}}
\newcommand{\econj}{\end{conjecture}}
\newcommand{\brem}{\begin{remark}}
\newcommand{\erem}{\end{remark}}

\newcommand{\ba}{\begin{array}}
\newcommand{\ea}{\end{array}}
\newcommand{\bea}{\begin{eqnarray}}
\newcommand{\eea}{\end{eqnarray}}
\newcommand{\bean}{\begin{eqnarray*}}
\newcommand{\eean}{\end{eqnarray*}}

\title{\huge On affine  orbifold nets associated with outer automorphisms\\}
\author{
{\sc Feng Xu}\footnote{Supported in part by NSF.}\\
Department of Mathematics\\
University of California at Riverside\\
Riverside, CA 92521\\
E-mail: {\tt xufeng@math.ucr.edu}}
\begin{document}
\date{}
\maketitle

\begin{abstract}
We construct  solitons in affine orbifold nets  associated with
outer automorphisms, and  we show that our construction gives all
the twisted representations of the fixed point subnet. This allows
us to settle a number of questions concerning such orbifold
constructions.
\end{abstract}

\newpage

\section{Introduction}

Let $\A$ be a completely rational conformal net (cf. \S
\ref{complete rationality} and def. \ref{abr} ). Let $\Gamma$ be a
finite group acting properly on $\A$ (cf. definition (\ref{p'})). It
is proved in   \cite{Xorb} that the fixed point subnet (the
orbifold) $\A^\Gamma$ is also completely rational, and by \cite{KLM}
$\A^\Gamma$ has finitely many irreducible representations which are
divided into two classes: the ones that are obtained from the
restrictions of a representation of $\A$ to $\A^\Gamma$ which are
called untwisted representations, and the ones which are twisted
(cf. definition after Th. \ref{orb}). It follows from
 Th. \ref{orb} that twisted representation of $\A^\Gamma$
always exists if  $\A^\Gamma\neq \A.$ The motivating question for
this paper is  to construct these   twisted representations of
$\A^\Gamma.$\par It turns out that all  representations of
$\A^\Gamma $ are closely related to the solitons of $\A$ (cf. \S3.3
and Prop. \ref{sg}). Solitons are representations of $\A_0$, the
restriction of $\A$ to the real line identified with a circle with
one point removed. Every representation of $\A$ restricts to a
soliton of $\A_0$, but not every soliton of $\A_0$ can be extended
to a representation of $\A$. The construction of soliton depends on
the net $\A$ and the action of $\Gamma.$  In this paper we consider
the orbifold net associated with $SU(n)_k$ with outer automorphism.
When $k=1$, this case is already covered by the general results of
\cite{DX} and in the framework of Vertex Algebras by \cite{BK}. Our
work will build upon the results of \cite{DX} and \cite{BK}.

The main difference between the case considered in this paper and
those of \cite{KLX} is that it is not an easy question to determine
the index of solitons, hence the strategy  of adding up all indices
to see if it agrees with the index formula in \cite{KLX} does not
work in the present case (cf. Cor. \ref{n3index} for the list of
indices when $n=3$).  Here we use a result of \cite{BEK} to count
the number of irreducible solitons. This allows us to show that the
list of known solitons constructed in section \ref{generalcase} is
in fact all the irreducible solitons and hence all irreducible
representations of the fixed point subnet (cf. Th. \ref{allirrep})
can be determined. We expect that this idea will work in other
cases.  Though our main results Th. \ref{allirrep} and Cor.
\ref{mtc} are expected from various partial results (cf. \cite{Ver})
, they have not appeared before, and we give some applications of
these results.\par

The rest of this paper is as follows: after preliminary sections on
nets and related concepts, we construct solitons using the ideas of
\cite{DX}. By a counting argument, we prove that all irreps of the
fixed point nets have been constructed in Th. \ref{allirrep}. In \S4
we consider applications of our main results in \S3 to the
properties of certain fusion matrices which have been studied from
different points of view in \cite{PZ} and \cite{BFS}  where these
properties are postulated motivated by considerations from boundary
conformal field theories. In Prop. \ref{n3fusion} we explicitly
determine these fusion matrices for the first non-trivial case when
$n=3$ using a result in \cite{CGS} about set of lines in Euclidean
space which mutually have the angles $\pi/3$ or $\pi/2.$ The result
agrees with formulas in \cite{PZ} and \cite{BFS}. As a corollary we
determine the set of indices of twisted solitons in Cor.
\ref{n3index}. It is an interesting question to extend these results
to $n\geq 4.$\par The author would like to thank Prof. V. G. Kac for
stimulating discussions and  providing references and useful
suggestions. The paper would not have been written without his help.
\section{Conformal nets on $S^1$}
\label{nets} In this section we introduce some basic concepts and
notations which will be used later. We refer the reader to \S2 of
\cite{KLX} for more details. \par
 By an interval of the circle we
mean an open connected non-empty subset $I$ of $S^1$ such that the
interior of its complement $I'$ is not empty. We denote by $\I$ the
family of all intervals of $S^1$.

A {\it net} $\A$ of von Neumann algebras on $S^1$ is a map
\[
I\in\I\to\A(I)\subset B(\H)
\]
from $\I$ to von Neumann algebras on a fixed Hilbert space $\H$ that
satisfies:
\begin{itemize}
\item[{\bf A.}] {\it Isotony}. If $I_{1}\subset I_{2}$ belong to
$\I$, then
\begin{equation*}
 \A(I_{1})\subset\A(I_{2}).
\end{equation*}
\end{itemize}
The net $\A$ is called {\it local} if it satisfies:
\begin{itemize}
\item[{\bf B.}] {\it Locality}. If $I_{1},I_{2}\in\I$ and $I_1\cap
I_2=\emptyset$ then
\begin{equation*}
 [\A(I_{1}),\A(I_{2})]=\{0\},
 \end{equation*}
where brackets denote the commutator.
\end{itemize}
The net $\A$ is called {\it M\"{o}bius covariant} if in addition
satisfies the following properties {\bf C,D,E,F}:
\begin{itemize}
\item[{\bf C.}] {\it M\"{o}bius covariance}.
There exists a strongly continuous unitary representation $U$ of the
M\"{o}bius group $\Mob$ (isomorphic to $PSU(1,1)$) on $\H$ such that
\begin{equation*}
 U(g)\A(I) U(g)^*\ =\ \A(gI),\quad g\in \Mob,\ I\in\I.
\end{equation*}
\end{itemize}
If $E\subset S^1$ is any region, we shall put
$\A(E)\equiv\bigvee_{E\supset I\in\I}\A(I)$ with $\A(E)=\mathbb C$
if $E$ has empty interior (the symbol $\bigvee$ denotes the von
Neumann algebra generated).  Note that the definition of $\A(E)$
remains the same if $E$ is an interval namely: if $\{I_n\}$ is an
increasing sequence of intervals and $\cup_n I_n = I$, then the
$\A(I_n)$'s generate $\A(I)$ (consider a sequence of elements
$g_n\in\Mob$ converging to the identity such that $g_n I\subset
I_n$).
\begin{itemize}
\item[{\bf D.}] {\it Positivity of the energy}.
The generator of the one-parameter rotation subgroup of $U$
(conformal Hamiltonian) is positive.
\item[{\bf E.}] {\it Existence of the vacuum}.  There exists a unit
$U$-invariant vector $\Omega\in\H$ (vacuum vector), and $\Omega$ is
cyclic for the von Neumann algebra $\bigvee_{I\in\I}\A(I)$.
\end{itemize}
By the Reeh-Schlieder theorem $\Omega$ is cyclic and separating for
every fixed $\A(I)$. The modular objects associated with
$(\A(I),\Omega)$ have a geometric meaning
\[
\Delta^{it}_I = U(\Lambda_I(2\pi t)),\qquad J_I = U(r_I)\ .
\]
Here $\Lambda_I$ is a canonical one-parameter subgroup of $\Mob$ and
$U(r_I)$ is a antiunitary acting geometrically on $\A$ as a
reflection $r_I$ on $S^1$.

This implies {\em Haag duality}:
\[
\A(I)'=\A(I'),\quad I\in\I\ ,
\]
where $I'$ is the interior of $S^1\setminus I$.

\begin{itemize}
\item[{\bf F.}] {\it Irreducibility}. $\bigvee_{I\in\I}\A(I)=B(\H)$.
Indeed $\A$ is irreducible iff $\Om$ is the unique $U$-invariant
vector (up to scalar multiples). Also  $\A$ is irreducible iff the
local von Neumann algebras $\A(I)$ are factors. In this case they
are III$_1$-factors in Connes classification of type III factors
(unless $\A(I)=\mathbb C$ for all $I$).
\end{itemize}
By a {\it conformal net} (or diffeomorphism covariant net) $\A$ we
shall mean a M\"{o}bius covariant net such that the following holds:
\begin{itemize}
\item[{\bf G.}] {\it Conformal covariance}. There exists a projective
unitary representation $U$ of $\Diff(S^1)$ on $\H$ extending the
unitary representation of $\Mob$ such that for all $I\in\I$ we have
\begin{gather*}
 U(g)\A(I) U(g)^*\ =\ \A(gI),\quad  g\in\Diff(S^1), \\
 U(g)xU(g)^*\ =\ x,\quad x\in\A(I),\ g\in\Diff(I'),
\end{gather*}
\end{itemize}
where $\Diff(S^1)$ denotes the group of smooth, positively oriented
diffeomorphism of $S^1$ and $\Diff(I)$ the subgroup of
diffeomorphisms $g$ such that $g(z)=z$ for all $z\in I'$.
\par
Let $G$ be a simply connected  compact Lie group. By Th. 3.2 of
\cite{FG}, the vacuum positive energy representation of the loop
group $LG$ (cf. \cite{PS}) at level $k$ gives rise to an irreducible
conformal net denoted by {\it ${\A}_{G_k}$}. By Th. 3.3 of
\cite{FG}, every irreducible positive energy representation of the
loop group $LG$ at level $k$ gives rise to  an irreducible covariant
representation of ${\A}_{G_k}$.
\subsection{Genus 0 $S,T$-matrices}
Next we will recall some of the results of \cite{R2}  and introduce
notations. \par Let $\{[\lambda], \lambda\in \L \}$ be a finite set
of all equivalence classes of irreducible, covariant, finite-index
representations of an irreducible local conformal net $\A$. We will
denote the conjugate of $[\lambda]$ by $[{\bar \lambda}]$ and
identity sector (corresponding to the vacuum representation) by
$[1]$ if no confusion arises, and let $N_{\lambda\mu}^\nu = \langle
[\lambda][\mu], [\nu]\rangle $. Here $\langle \mu,\nu\rangle$
denotes the dimension of the space of intertwiners from $\mu$ to
$\nu$ (denoted by $\text {\rm Hom}(\mu,\nu)$).  We will denote by
$\{T_e\}$ a basis of isometries in $\text {\rm
Hom}(\nu,\lambda\mu)$. The univalence of $\lambda$ and the
statistical dimension of (cf. \S2  of \cite{GL1}) will be denoted by
$\omega_{\lambda}$ and $d{(\lambda)}$ (or $d_{\lambda})$)
respectively. \par Let $\phi_\lambda$ be the unique minimal left
inverse of $\lambda$, define:
\begin{equation}\label{Ymatrix}
Y_{\lambda\mu}:= d(\lambda)  d(\mu) \phi_\mu (\epsilon (\mu,
\lambda)^* \epsilon (\lambda, \mu)^*),
\end{equation}
where $\epsilon (\mu, \lambda)$ is the unitary braiding operator
 (cf. \cite{GL1} ). \par
We list two properties of $Y_{\lambda \mu}$ (cf. (5.13), (5.14) of
\cite{R2}) which will be used in the following:
\begin{lemma}\label{Yprop}
\begin{equation*}
Y_{\lambda\mu} = Y_{\mu\lambda}  = Y_{\lambda\bar \mu}^* = Y_{\bar
\lambda \bar \mu}.
\end{equation*}
\begin{equation*}
Y_{\lambda\mu}  = \sum_k N_{\lambda\mu}^\nu
\frac{\omega_\lambda\omega_\mu} {\omega_\nu} d(\nu) .
\end{equation*}
\end{lemma}
We note that one may take the second equation in the above lemma as
the definition of $Y_{\lambda\mu}$.\par Define $a := \sum_i
d_{\rho_i}^2 \omega_{\rho_i}^{-1}$. If the matrix $(Y_{\mu\nu})$ is
invertible, by Proposition on P.351 of \cite{R2} $a$ satisfies
$|a|^2 = \sum_\lambda d(\lambda)^2$.
\begin{definition}\label{c0}
Let $a= |a| \exp(-2\pi i \frac{c_0}{8})$ where  $c_0\in {\mathbb R}$
and $c_0$ is well defined ${\rm mod} \ 8\mathbb Z$.
\end{definition}
Define matrices
\begin{equation}\label{Smatrix}
S:= |a|^{-1} Y, T:=  C {\rm Diag}(\omega_{\lambda})
\end{equation}
where \[C:= \label{dims} \exp(-2\pi i \frac{c_0}{24}).\] Then these
matrices satisfy (cf. \cite{R2}):
\begin{lemma}\label{Sprop}
\begin{align*}
SS^{\dag} & = TT^{\dag} ={\rm id},  \\
STS &= T^{-1}ST^{-1},  \\
S^2 & =\hat{C},\\
 T\hat{C} & =\hat{C}T=T,
\end{align*}

where $\hat{C}_{\lambda\mu} = \delta_{\lambda\bar \mu}$ is the
conjugation matrix.
\end{lemma}
Moreover
\begin{equation}\label{Verlinde}
N_{\lambda\mu}^\nu = \sum_\delta \frac{S_{\lambda\delta}
S_{\mu\delta} S_{\nu\delta}^*}{S_{1\delta}}. \
\end{equation}
is known as Verlinde formula. \par We will refer the $S,T$ matrices
as defined above  as  {\bf genus 0 modular matrices of ${\A}$} since
they are constructed from the fusion rules, monodromies and minimal
indices which can be thought as  genus 0 {\bf chiral data}
associated to a Conformal Field Theory. \par The commutative algebra
generated by $\lambda$'s with structure constants
$N_{\lambda\mu}^\nu$ is called {\bf fusion algebra} of $\A$. If $Y$
is invertible, it follows from Lemma \ref{Sprop}, (\ref{Verlinde})
that any nontrivial irreducible representation of the fusion algebra
is of the form $\lambda\rightarrow \frac{S_{\lambda\mu}}{S_{1\mu}}$
for some $\mu$.
\par
\subsection{The orbifolds}
Let ${\A}$ be an irreducible conformal net on a Hilbert space ${\H}$
and let $\Gamma$ be a finite group. Let $V:\Gamma\rightarrow
U({\H})$ be a  unitary representation of $\Gamma$ on ${\H}$. If
$V:\Gamma\rightarrow U({\H})$ is not faithful, we set $\Gamma':=
\Gamma /{\rm ker} V$.
\begin{definition} \label{p'}
We say that $\Gamma$ acts properly on ${\A}$ if the following
conditions are satisfied:\par (1) For each fixed interval $I$ and
each $g\in \Gamma$, $\alpha_g (a):=V(g)aV(g^*) \in {\A}(I), \forall
a\in {\A}(I)$; \par (2) For each  $g\in \Gamma$, $V(g)\Omega =
\Omega, \forall g\in \Gamma$.\par \label{Definition 2.1}
\end{definition}
We note that if $\Gamma$ acts properly, then $V(g)$, $g\in\Gamma$
commutes with the unitary representation $U$ of $\Mob$. \par Define
$\B(I):= \{ a\in \A(I) | \alpha_g (a)=a, \forall g\in \Gamma \}$ and
${\A}^\Gamma(I):={\B}(I)P_0$ on ${\H}_0$ where $\H_0:=\{ x\in \H|
V(g)x=x, \forall g\in \Gamma \}$ and $P_0$ is the projection from
$\H$ to $\H_0.$  Then $U$ restricts to an  unitary representation
(still denoted by $U$) of $\Mob$ on ${\H}_0$. The following is
proved in \cite{Xorb}: \bprop\label{Prop.2.1}  The map $I\in
{\I}\rightarrow {\A}^{\Gamma}(I)$ on $ {\H}_0$ together with the
unitary representation (still denoted by $U$) of \Mob\ on ${\H}_0$
is an irreducible M\"{o}bius covariant net.  \eprop The irreducible
M\"{o}bius covariant net in Prop. \ref{Prop.2.1} will be denoted by
${\A}^\Gamma$ and will be called the {\it orbifold of ${\A}$} with
respect to $\Gamma$. We note that by definition ${\A}^\Gamma=
{\A}^{\Gamma'}$.
\par
\subsection{Complete rationality }
\label{complete rationality} We first recall some definitions from
\cite{KLM} . Recall that   ${\I}$ denotes the set of intervals of
$S^1$. Let $I_1, I_2\in {\I}$. We say that $I_1, I_2$ are disjoint
if $\bar I_1\cap \bar I_2=\emptyset$, where $\bar I$ is the closure
of $I$ in $S^1$. When $I_1, I_2$ are disjoint, $I_1\cup I_2$ is
called a 1-disconnected interval in \cite{Xjw}. Denote by ${\I}_2$
the set of unions of disjoint 2 elements in ${\I}$. Let ${\A}$ be an
irreducible M\"{o}bius covariant net as in \S2.1. For $E=I_1\cup
I_2\in{\I}_2$, let $I_3\cup I_4$ be the interior of the complement
of $I_1\cup I_2$ in $S^1$ where $I_3, I_4$ are disjoint intervals.
Let
$$
{\A}(E):= A(I_1)\bigvee A(I_2), \quad \hat {\A}(E):= (A(I_3)\bigvee
A(I_4))'.
$$ Note that ${\A}(E) \subset \hat {\A}(E)$.
Recall that a net ${\A}$ is {\it split} if ${\A}(I_1)\bigvee
{\A}(I_2)$ is naturally isomorphic to the tensor product of von
Neumann algebras ${\A}(I_1)\otimes {\A}(I_2)$ for any disjoint
intervals $I_1, I_2\in {\I}$. ${\A}$ is {\it strongly additive} if
${\A}(I_1)\bigvee {\A}(I_2)= {\A}(I)$ where $I_1\cup I_2$ is
obtained by removing an interior point from $I$. \bdefin\label{abr}
\cite{KLM} ${\A}$ is said to be completely  rational if ${\A}$ is
split, strongly additive, and the index $[\hat {\A}(E): {\A}(E)]$ is
finite for some $E\in {\I}_2$ . The value of the index $[\hat
{\A}(E): {\A}(E)]$ (it is independent of $E$ by Prop. 5 of
\cite{KLM}) is denoted by $\mu_{{\A}}$ and is called the $\mu$-index
of ${\A}$. If the index $[\hat {\A}(E): {\A}(E)]$ is infinity for
some $E\in {\I}_2$, we define the $\mu$-index of ${\A}$ to be
infinity. \label{Definition 2.2} \edefin

Note that by \cite{LX} every irreducible, split, local conformal net
with finite $\mu$-index is automatically strongly additive. The
following theorem is proved in \cite{Xorb}: \btheor\label{orb} Let
${\A}$ be an irreducible M\"{o}bius covariant net and let $\Gamma$
be a finite group acting properly on ${\A}$. Suppose that ${\A}$ is
completely rational. Then:\par (1): ${\A}^\Gamma$ is completely
rational or $\mu$-rational and $\mu_{{\A}^\Gamma}= |\Gamma'|^2
\mu_{{\A}}$; \par (2): There are only a finite number of irreducible
covariant representations of ${\A}^\Gamma$ (up to unitary
equivalence), and they give rise to a unitary modular category as
defined in II.5 of \cite{Tu} by the construction as given in \S1.7
of \cite{X3m}. \label{Th.2.6} \etheor Suppose that ${\A}$ and
$\Gamma$ satisfy the assumptions of Th. \ref{orb}. Then
${\A}^\Gamma$ has only finitely number of irreducible
representations $\dot\lambda$ and
$$
\sum_{\dot\lambda}d(\dot\lambda)^2 = \mu_{{\A}^\Gamma}= |\Gamma'|^2
\mu_{{\A}} .
$$

\subsection{Restriction to the real line: Solitons}
Denote by $\I_0$ the set of open, connected, non-empty, proper
subsets of $\mathbb R$, thus $I\in\I_0$ iff $I$ is an open interval
or half-line (by an interval of $\mathbb R$ we shall always mean a
non-empty open bounded interval of $\mathbb R$).

Given a net $\A$ on $S^1$ we shall denote by $\A_0$ its restriction
to $\mathbb R = S^1\setminus\{-1\}$. Thus $\A_0$ is an isotone map
on $\I_0$, that we call a \emph{net on $\mathbb R$}. In this paper
we denote by $J_0:=(0,\infty)\subset \mathbb R$.

A representation $\pi$ of $\A_0$ on a Hilbert space $\H$ is a map
$I\in\I_0\mapsto\pi_I$ that associates to each $I\in\I_0$ a normal
representation of $\A(I)$ on $B(\H)$ such that
\[
\pi_{\tilde I}\res\A(I)=\pi_I,\quad I\subset\tilde I, \quad I,\tilde
I\in\I_0\ .
\]
A representation $\pi$ of $\A_0$ is also called a \emph{soliton}. As
$\A_0$ satisfies half-line duality, namely
$$
\A_0(-\infty,a)'= \A_0(a,\infty), a\in \mathbb R,
$$
by the usual DHR argument \cite{DHR} $\pi$ is unitarily equivalent
to a representation $\rho$ which  acts identically on
$\A_0(-\infty,0)$, thus $\rho$ restricts to an endomorphism of
$\A(J_0)= \A_0(0,\infty)$. $\rho$ is said to be localized on $J_0$
and we also refer to $\rho$ as soliton endomorphism.

Clearly a representation $\pi$ of $\A$ restricts to a soliton
$\pi_0$ of $\A_0$. But a representation $\pi_0$ of $\A_0$ does not
necessarily extend to a representation of $\A$.

If $\A$ is strongly additive, and a representation $\pi_0$ of $\A_0$
extends to a DHR representation of $\A$, then it is easy to see that
such an extension is unique, and in this case we will use the same
notation $\pi_0$ to denote the corresponding DHR  representation of
$\A$.

\subsection{Induction and restriction}

Let $\A$ be a M\"obius covariant net and $\B$ a subnet. Given  a
bounded interval $I_0\in\I_0$ we fix canonical endomorphism
$\gamma_{I_0}$ associated with $\B(I_0)\subset\A(I_0)$. Then we can
choose for each $I\in\I_0$ with $I\supset I_0$ a canonical
endomorphism $\gamma_{I}$ of $\A(I)$ into $\B(I)$ in such a way that
$\gamma_{I}\res\A(I_0)=\gamma_{I_0}$ and $\l_{I_1}$ is the identity
on $\B(I_1)$ if $I_1\in\I_0$ is disjoint from $I_0$, where
$\l_{I}\equiv\gamma_{I}\res\B(I)$.

We then have an endomorphism $\gamma$ of the $C^*$-algebra
$\mathfrak A\equiv\overline{\cup_{I}\A(I)}$ ($I$ bounded interval of
$\mathbb R$).

Given a DHR endomorphism $\r$ of $\B$ localized in $I_0$, the
$\a$-induction $\a_{\r}$ of $\r$ is the endomorphism of $\mathfrak
A$ given by
\[
\a_{\r}\equiv \gamma^{-1}\cdot\Ad\e(\r,\l)\cdot\r\cdot\gamma\ ,
\]
where $\e$ denotes the right braiding unitary symmetry (there is
another choice for $\a$ associated with the left braiding).
$\a_{\r}$ is localized in a right half-line containing $I_0$, namely
$\a_\r$ is the identity on $\A(I)$ if $I$ is a bounded interval
contained in the left complement of $I_0$ in $\mathbb R$. Up to
unitarily equivalence, $\a_\r$ is localizable in any right half-line
thus $\a_\r$ is normal on left half-lines, that is to say, for every
$a\in\mathbb R$, $\a_\r$ is normal on the $C^*$-algebra $\mathfrak
A(-\infty,a)\equiv\overline{\cup_{I\subset (-\infty,a)}\A(I)}$ ($I$
bounded interval of $\mathbb R$), namely $\a_\r\res\mathfrak
A(-\infty,a)$ extends to a normal morphism of $\A(-\infty,a)$. We
have the following Prop. 3.1 of \cite{LX}:
\begin{proposition}\label{sg}
$\a_\r$ is a soliton endomorphism of $\A_0$.
\end{proposition}
\subsection {Loop
groups of type A}\label{typea} We denote $LSU(n)$ the group of
smooth maps $f: S^1 \mapsto SU(n)$ under pointwise multiplication.
The diffeomorphism group of the circle $\text{\rm Diff} S^1 $ is
naturally a subgroup of $\text{\rm Aut}(LSU(n))$ with the action
given by reparametrization. In particular the group of rotations
$\text{\rm Rot}S^1 \simeq U(1)$ acts on $LSU(n)$. The Lie algebra of
$LSU(n)$, denoted by $Lsu(n),$ consists of smooth maps from $S^1$ to
$su(n).$  We will denote elements of $Lsu(n)$ by its Fourier series
$g(z)= \sum_{n} g_n z^n,$ and $L^0su(n)$ the subspace of $Lsu(n)$
which are polynomials in $z=\exp(2\pi i \theta), 0\leq \theta\leq
1.$  We will be interested in the projective unitary representation
$\pi : LSU(n) \rightarrow U(H)$ that are both irreducible and have
positive energy. This means that $\pi $ should extend to
$LSU(n)\ltimes \text{\rm Rot}\ S^1$ so that $H=\oplus _{n\geq 0}
H(n)$, where the $H(n)$ are the eigenspace for the action of
$\text{\rm Rot}S^1$, i.e., $r_\theta \xi = \exp(i n \theta)$ for
$\theta \in H(n)$ and $\text{\rm dim}\ H(n) < \infty $ with $H(0)
\neq 0$. By I. 7 of \cite{Wa} the space of finite energy vectors are
$C^\infty$ vectors for the action of $L^0su(n)$, and by I.9 of
\cite{Wa} $H$ remains irreducible when restricting to subgroups
generated by $\exp(i X), X=X^*\in L^0su(n).$ We will use $\L SU(n)$
to denote the central extension of $LSU(n)$ by $S^1$ as constructed
in Chapter 4 of \cite{PS}.

It follows from \cite{PS} and \cite{Kac} that for fixed level $k$
which is a positive integer, there are only finite number of such
irreducible representations indexed by the finite set
$$
 P_{++}^{k}
= \bigg \{ \lambda \in P \mid \lambda = \sum _{i=1, \cdots , n-1}
\lambda _i \Lambda _i , \lambda _i \geq 0\, , \sum _{i=1, \cdots ,
n-1} \lambda _i \leq k \bigg \}
$$
where $P$ is the weight lattice of $SU(n)$ and $\Lambda _i$ are the
fundamental weights.

We will use $\Lambda_0$ or simply $1$  to denote the trivial
representation of $SU(n)$. For $\lambda , \mu , \nu \in P_{++}^{k}$,
define $N_{\lambda \mu}^\nu  = \sum _{\delta \in P_{++}^{k}
}S_\lambda ^{(\delta)} S_\mu ^{(\delta)} S_\nu
^{(\delta*)}/S_{\Lambda_0}^{(\delta)}$ where $S_\lambda ^{(\delta)}$
is given by the Kac-Peterson formula (cf. equation (\ref{kacp})
below for an equivalent formula):
$$
S_\lambda ^{(\delta)} = c \sum _{w\in S_n} \varepsilon _w \exp
(iw(\delta) \cdot \lambda 2 \pi /n)
$$
where $\varepsilon _w = \text{\rm det}(w)$ and $c$ is a
normalization constant fixed by the requirement that
$S_\mu^{(\delta)}$ is an orthonormal system. It is shown in
\cite{Kac} P. 288 that $N_{\lambda \mu}^\nu $ are non-negative
integers. Moreover, define $ Gr(C_k)$ to be the ring whose basis are
elements of $ P_{++}^{k}$ with structure constants $N_{\lambda
\mu}^\nu $.
  The natural involution $*$ on $ P_{++}^{k}$ is
defined by $\lambda \mapsto \bar{\lambda}  =$ the conjugate of
$\lambda $ as representation of $SU(n)$.\par

We shall also denote $S_{\Lambda _0}^{(\Lambda)}$ by $S_1^{(\Lambda
)}$. Define $d_\lambda = \frac {S_1^{(\lambda )}}{S_1^{(\Lambda
_0)}}$. We shall call $(S_\nu ^{(\delta )})$ the $S$-matrix of
$LSU(n)$ at level $k$. \par

The irreducible positive energy representations of $ L SU(n)$ at
level $k$ give rise to an irreducible conformal net $\A$ (cf.
\cite{KLX}) and its covariant representations.

We will use $\lambda=(\lambda_1,...\lambda_{n-1})$ to denote
irreducible representations of $\A$ and also the corresponding
endomorphism of $M=\A(I).$ Recall from \cite{KLX} that $\A(I)$ is
generated as von Neumann algebra by $\pi_{k\Lambda_0}(f), \forall
f\in LSU(n), f\res I' = e$ where $e$ denotes the identity element of
$SU(n).$

All the sectors $[\lambda]$ with $\lambda$ irreducible generate the
fusion ring of $\A.$
\par

The following form of Kac-Peterson formula for $S$ matrix will be
used later:
\begin{equation}\label{kacp}
\frac{S_{\lambda\mu}}{S_{1\mu}}= \exp(\frac{t(\mu+\rho)}{n(k+n)})
{ch}_{\lambda'} (x_1,...,x_{n-1}, 1)
\end{equation}
Where ${ch}_{\lambda'}$ is the character associated with finite
irreducible  representation of $SU(n)$ labeled by $\lambda,$ and
$x_i=\exp(-2\pi i\frac{\mu_i'}{k+n}), \mu_i'=\sum_ {i\leq j\leq n-1}
(\mu_j+1), 1\leq i\leq n-1, t(\lambda)= \sum_{1\leq i\leq n-1}
i\lambda_i.$
\par

The following  result is proved in \cite{Wa} (See Corollary 1 of
Chapter V in \cite{Wa}).

\btheor\label{wass}  Each $\lambda \in  P_{++}^{(k)}$ has finite
index with index value $d_\lambda ^2$.  The fusion ring generated by
all $\lambda \in P_{++}^{(k)}$ is isomorphic to $ Gr(C_k)$. \etheor
In the case of $SU(2)_k,$ we will label irreducible representations
by a half integer called spin $0\leq i\leq k/2.$ Here are some
examples of fusion rules:
\begin{equation}\label{su2}
\frac{1}{2}\times i = (i-\frac{1}{2})\oplus (i+\frac{1}{2}), 0\leq
i\leq \frac{k-1}{2};   1\times i = (i-1)\oplus i\oplus (i+1),0\leq
i\leq \frac{k-2}{2}.
\end{equation}
\subsection{Twisted loop group}
Let $\tau$ be the order two outer automorphism of $SU(n)$ given by
$\t(A)=\bar{A}, \forall A\in SU(n)$ where $\bar{A}$ is the complex
conjugate of $A.$ On the Lie algebra $su(n)$ $\tau$ is given by
$\tau(X)= -\bar X, \forall X\in su(n)$ (We identify $su(n)$ with
$n\times n$ Hermitian matrices). It is convenient in this paper to
think of twisted loop group $L_\t SU(n)$ as a subgroup of $LSU(n).$
We make the following definition: \bdefin\label{twistedloop} $L_\t
SU(n):= \{ f\in LSU(n), f(\theta + \frac{1}{2})= \tau (f(\theta)),
0\leq\theta\leq \frac{1}{2}\};$ $L_\t su(n):= \{ f\in Lsu(n),
f(\theta + \frac{1}{2})= \tau (f(\theta)), 0\leq\theta\leq
\frac{1}{2}\}.$ \edefin
Let $su(n)= so(n)\oplus g_1$ where $g_1$ is the eigenspace for $\t$
with eigenvalue $-1.$ Note that $g_1$ is an irreducible
representation of $so(n)$ under the adjoint action (cf. \S8 of
\cite{Kac}). Assume that $t=t_0\oplus t_1$ where $t\subset su(n)$ is
the subset of diagonal matrices. Note that both $t_0,t_1$ are
nontrivial subspaces.

\par The
projective irreducible representations of $L_\t SU(n)$ are similar
to that of $LSU(n).$ For each fixed positive integer $k$ we again
have finitely many projective unitary irreducible representations of
$L_\t SU(n)\ltimes S^1.$ We shall use $\L_\t SU(n)$ (resp. $\L_\t
su(n)$) to denote the central extensions of $L_\t SU(n)$ (resp.
$L_\t su(n)$ ) so that the projective unitary irreducible
representations of $L_\t SU(n)\ltimes S^1$ at level $k$ are
irreducible representations of $\L_\t SU(n)\ltimes S^1.$ These
representations correspond to irreducible integrable highest weight
modules of $\L_\t su(n)$  at level $k.$ We refer the reader to
Chapter II of \cite{Ver} for more details. The following simple
observation follows from Chapter 10 of \cite{Kac}:

\blemma\label{tirrep} The number of  irreducible representations of
$\L_\t SU(n)$ at level $k$ is equal to the number of irreducible
representations $\lambda$ of $\L SU(n)$ at level $k$ such that
$\lambda=\bar{\lambda}.$\elemma

\blemma\label{cartan} Let $\B$ be the subnet of $\A_{SU(n)_k}$ such
that $\B(I)$ is generated as a von Neumann algebra by
$\pi_{k\Lambda_0} (g\exp(i X) g^*), \forall  g\in Spin(n), X=X^*\in
Lt, X(z)=0, \forall z\in I'$. Then $\B=\A$.\elemma \proof Note that
$\B$ is indeed a subnet of $\A$ since the modular group associated
with $\A(I)$ preserves $\B(I).$ To show that $\B(I)=\A(I)$ it is
enough to show that $\bigvee_{I} \B(I)=\bigvee_I\A(I).$ Let $Y=gyg^*
z^n + gy^*g^* z^{-n}\in Lsu(n)$ with $y\in t$ and $g\in Spin(n).$
Then $\pi_{k\Lambda_0}(\exp{iY}) \in \bigvee_{I} \B(I).$ Since the
adjoint actions of $Spin(n)$ on $g_0, g_1$ are irreducible, and
$t_0,t_1$ are nontrivial subspaces, it follows that $g_0$ (resp.
$g_1$) is the linear span of $g X g^*, \forall g \in Spin(n), X\neq
0 \in t_0$ (resp. $g X g^*, \forall g \in Spin(n),  X\neq 0 \in
t_1$). By Trotter's product formulas (cf. P. 295 of \cite{RS}) we
conclude that $\pi_{k\Lambda_0}(\exp{iY}) \in \bigvee_{I} \B(I)$ for
any $Y=Y^*\in L^{0}su(n).$ Since $\pi_{k\Lambda_0}$ is irreducible
as representation of the group generated by $\exp(iY), \forall
Y=Y^*\in L^{0}su(n),$ the lemma is proved.\endproof

\section{Affine orbifold nets associated with outer automorphisms}
Let $k$ be a positive integer (level). Unless otherwise stated we
will write $\A= \A_{SU(n)_k}.$ By identifying $\R^{2n}=(x,y)
\rightarrow x+iy\in {\Bbb C}^{n}$ where $x,y$ are column vectors
with $n$ real entries, we have the natural inclusion $SU(n)_k\subset
Spin(n)_k.$ Define $J:=(Id_n, -Id_n)\in SO(2n)$ and lift it to
$Spin(2n).$ Note that for $A\in SU(N), JAJ=\bar{A},$ and $J. \Omega=
\Omega$ where $\Omega$ is the vacuum vector for vacuum
representations of $\L Spin(2n)_k.$ It follows that $Ad J$ generates
a proper $\Z_2$ action of $\A_{SU(n)_k}.$ This is the action
corresponding to the outer automorphism $\tau$ of $SU(n).$ Suppose
that the vacuum representation of $\A$ decompose as $1\oplus \sigma$
as representations of $\A^{\Z_2}$ where $1$ stands for the vacuum
representation. Motivated by Lemma 8.3 of \cite{LX}, we make the
following definition: \bdefin\label{ttwisted} Let $\rho$ be an
irreducible soliton of $\A$ which restricts to a DHR representation
of $\A^{\Bbb Z_2},$ and let $\rho_+$ be an irreducible component of $\rho\res
\A^{\Bbb Z2}.$ $\rho$ is called $\t$-twisted if $ \e(\rho_+,\sigma)\e(\sigma,
\rho_+)= -1 $ where $\e(.,.)$ is the braiding operator (cf.
\cite{GL1}). \edefin

\subsection{Constructions of solitons }

\bdefin\label{lt}
$$L_\R SU(n):= \{
f\in LSU(n)| f(0)=f(1)= e, f^{(n)}(0)= f^{(n)}(1)=0, \forall n\geq 1
\};
$$
$$L_\R T_0:= \{ \exp(i f(e^{2\pi i\theta})),
f=f^*\in Lt| f^{(n)}(0)= f^{(n)}(1)=0, \forall n\geq 0 \}; $$
$$L_I
T_0:= \{ \exp(i f), f=f^*\in L t, I\subset \R, f\res I'= 0 \}. $$
\edefin The following is a special case  of covering homomorphism in
Prop. 4.6 of \cite{DX}: \bdefin\label{cover} $\phi:L_\R
T_0\rightarrow L_\t T$ is a homomorphism defined by $\phi(g)(\theta)=
g(2\theta), \phi(g)(\theta+\frac{1}{2})=\t(\phi(g)(\theta)),0\leq
\theta\leq 1/2.$ $\phi$ lifts to a homorphism from central
extensions of $\L_R T_0$ to central extensions of $\L_\t T$  and by
abuse of notations we will use $\phi$ to denote the lift. \edefin
\subsubsection{level 1 case}
When the level $k=1$, we have conformal inclusion $Spin (n)_2\subset
SU(n)_1,$ and by Lemma 5.1 of \cite{Xorb} we have $\A^{\Z_2}=
\A_{Spin(n)_2}.$ By 4.3 of \cite{Xorb} $\A_{Spin(n)_2}$ is
completely rational, and its irreducible representations are in one
to one correspondence to irreducible representations of $\L
Spin(n)_2.$ \par Denote by $T$ the subgroup of diagonal matrices in
$SU(n)$ and $\L T$  the subgroup of $\L SU(n)_1.$  The conformal net
$\A_T$ associated with $\L T$ is the same as $\A_{SU(n)_1}$ by Page
28 of \cite{Xcos}, and is a special case of conformal net associated
with lattices as defined in definition 3.7 of \cite{DX}. Hence the
$\t$-twisted irreducible solitons of $\A$ are given by section 4.1
of \cite{DX}.  These  $\t$-twisted solitons  are finitely direct sum
of irreducible representations of $\L Spin(n)_2$ (In particular they
furnish a representation of $Spin(n)$),  and they are in one to one
correspondence to irreducible representations of $\L_\t SU(n)$ at
level $1$ (cf. \cite{BK},\cite{FJ} and references therein). We
summarize these results in the following:

\blemma\label{level1} The list $\{ \pi_\t \}$ of irreducible
$\t$-twisted solitons of $\A$ is given as follows: they are in one
to one correspondence with the list of irreducible representations
$\{ \pi \}$ of $\L_\t SU(n)$ at level $1$, and we have
\begin{equation}\label{key}
\pi_\t (gfg^*) = \pi(g)\pi (\phi(f))\pi(g^*), \forall f\in \L_I
T_0\subset \L_\R T_0, \forall g\in Spin(n).
\end{equation}
\elemma

\blemma\label{dense} Let $X=gxg^*z^n+gx^*g^*z^{-n}\in L_\sigma
su(n)$ with $g\in Spin(n), x\in t.$ Then there exists a sequence of
$g_m\in \L_\R SU(n)$ such that
$$
\pi(\exp({iX}))= s-\lim_{m\rightarrow \infty} \pi(\phi(g_m))
$$
where $\pi $ is a direct sum of finitely many irreducible
representations of $\L_\t SU(n)_1.$ The same is true if $\pi$ is
replaced by $\pi \otimes \pi...\otimes \pi$ where there are $k$
tensor product of $\pi.$ \elemma \proof  The proof is similar to the
proof of Prop. 4.10 of \cite{DX}. By Prop. 1.2.3 in Chapter 4 of
\cite{TL} we can choose a sequence  $x_m(\theta)$ of smooth complex
valued functions on $[0,1]$ such that $x_m(\theta + \frac{1}{2})=
x_m(\theta), x_m^{(n)}(0)= x_m^{(n)}(1)=0, \forall n\geq 0, ||x_m
-1||_{\frac{1}{2}}\leq 1/m,  m\geq 1$ (cf. \S1.2 of \cite{TL} for
the definition of norm $||.||_{\frac{1}{2}}$). It follows by
definition that $\exp (i x_m X)= \phi(\exp (i f_m))$ where
$f_m(\theta)= x_m(\theta/2)X(\theta/2), 0\leq \theta \leq 1.$ Note
that $\pi$ is a  direct sum of finitely many irreducible
representations of $\L Spin(n)_2,$ by Prop. 1.3.2 in Chapter 4 of
\cite{TL} we have that $\pi(\phi(\exp(i f_m)))\rightarrow
\pi(\exp(iX))$ strongly. When $\pi$ is replaced by $\pi \otimes
\pi...\otimes \pi$ where there are $k$ tensor product of $\pi,$ the
same argument as in Prop. 1.3.2 in Chapter 4 of \cite{TL} works,
provided that  one replaces the generator $L_0$ by $\sum_ {1\leq
i\leq k} id\otimes ... \otimes L_0...\otimes id $ where in the
summation $ L_0$ appears in the $i$-th tensor.

\endproof

\subsection{General level case}\label{generalcase}
Let $\pi_1$ be a direct sum of all level $1$ irreducible
representations of $\L_\tau SU(n)$ and let $\pi$ be $k$ tensor
products of $\pi_1.$ Note that $\pi$ is a representation of $\L_\tau
SU(N)$ at level $k$ with positive energy. By Lemma \ref{level1}
$\pi$ gives a soliton of $\A_{SU(n)_1}^{\otimes k}\supset
\A_{SU(n)_k}$, and by restriction, a soliton $\pi_\tau$ of
$\A_{SU(n)_k}.$ We have
$$
\pi_\tau (gfg^*)= \pi(g)\pi(\phi(f))\pi(g^*), \forall f\in \L_I T_0,
I\subset \R, \forall g\in Spin(n),
$$
where $\phi$ is defined in definition \ref{cover}, and we have
identified $f$ with its image in $\A(I).$\par Note that
$\A_{SU(n)_1}^{\otimes k}= \A_T^{\otimes k}$ is a net associated
with $\L (T\times T...\times T)$ where there are $k$ products, and
is a net associated with lattice as in Definition 3.7 of \cite{DX}.
By Prop. 4.8 of \cite{DX} $\pi$ restricts to a DHR representation of
of $(\A_{SU(n)_1}^{\otimes k})^{\Z_2}$ where $\Z_2$ is generated by
$\t\otimes \t...\otimes\t,$ hence $\pi_\tau$ restricts to a DHR
representation of of $\A_{SU(n)_k}^{\Z_2}.$

\bprop\label{dense2}
$$
\bigvee_{I\subset \R} \pi_\tau(\A_I))= \pi (\L_\tau SU(n))''.
$$
\eprop \proof Since by Lemma \ref{cartan} $\A(I)$ is generated by $g
f g^*, \forall f\in \L_I T_0, \forall g\in Spin(n),$ by definition
we have
$$
\bigvee_{I\subset \R} \pi_\tau(\A_I))\subset  \pi (\L_\tau SU(n))''.
$$
Since $L_\tau SU(n)$ is connected (cf. Lemma 4.2 of \cite{Ver}), it
is sufficient to check that for $X\in L_\tau su(n)$,
$\pi(\exp(iX))\in \bigvee_{I\subset \R} \pi_\tau(\A_I).$ As in the
proof of Lemma \ref{cartan}, by Trotter's product formula (cf. Page
295 of \cite{RS}) and irreducibility of the actions of $Spin(n)$ on
$so(n)$ and $g_1,$ it is sufficient to check that for
$X=gxz^ng^*+gx^*z^{-n}g^*\in L_\sigma su(n)$ with $g\in Spin(n),
x\in t,$ we have $\pi(\exp(iX))\in \bigvee_{I\subset \R}
\pi_\tau(\A_I),$ and this follows from Lemma \ref{dense}.
\endproof

\bcoro\label{irrep} (1) Each irreducible representation $\rho$ of
$\L_\tau SU(n)$ at level $k$ gives an irreducible  soliton
$\rho_\tau$ of $\A_{SU(n)_k}$ such that $\rho_1\simeq \rho_2$ as
representations of $\L_\tau SU(n)$ at level $k$ iff
$\rho_{1\tau}\simeq \rho_{2\tau}$ as  solitons of
$\A_{SU(n)_k}$;\par (2) $\rho_\tau \tau\simeq \rho_\tau$ as solitons
of $\A_{SU(n)_k}$, and $\rho_\t$ restricts to a DHR representation
of $\A_{SU(n)_k}^{\Z_2}.$\ecoro \proof By \cite{Kac} all irreducible
representations of $\L_\tau SU(n)$ at level $k$ appear in $\pi,$ and
(1) follows from Prop. \ref{dense2}. For (2), first we note that
$\rho_\t$ comes from an irreducible component of $\pi_\t$ in Prop.
\ref{dense2}. Since $\pi_\t$ restricts to a DHR representation of
$\A_{SU(n)_k}^{\Z_2},$ it follows that $\rho_\t$ restricts to a DHR
representation of $\A_{SU(n)_k}^{\Z_2}.$ By construction
$$
\rho_\tau (\tau (gfg^*))= \rho(\phi(\tau(gfg^*)))= \rho( R_{1/2}
\phi(gfg^*)), \forall g\in Spin(n), f\in \L_I T_0
$$
where $R_{1/2} (g)(\theta) = g(\theta +\frac{1}{2}).$ Since
rotations are implemented on $\rho$ (cf. Page 246 in \cite{Kac} for
a formula for the infinitesimal generator of rotations), (2)
follows. \endproof
 Let $\rho_\tau$ be an irreducible soliton of $\A$
as given by Cor. \ref{irrep}. By (2) of Cor. \ref{allirrep} and Cor.
4.9 of \cite{KLM} $\rho_\tau$ decomposes into direct sum of two
distinct irreducible representations $\rho_{+}, \rho_-$ of
$\A^{\Z_2}. $ By the same argument as Prop. 4.17 of \cite{DX} it
follows that $\e(\rho_+,\sigma) \e(\sigma,\rho_+)= -1,$ and so we
have the following lemma: \blemma\label{tt} The irreducible solitons
as given in Cor. \ref{irrep} are $\t$-twisted as in definition
\ref{ttwisted}.\elemma
\subsection{Counting of all irreducible $\t$-twisted solitons}
Now we apply induction and restriction for general orbifolds in \S4
of \cite{KLM} to $\A^{\Z_2}\subset \A.$ Recall  that the vacuum
representation of $\A$ decompose as $1\oplus \sigma$ as
representations of $\A^{\Z_2}$ where $1$ stands for the vacuum
representation. Let $\rho_\tau$ be $\tau$-twisted irreducible
soliton of $\A$ as given by Cor. \ref{irrep}. By (2) of Cor.
\ref{allirrep} and Cor. 4.9 of \cite{KLM} $\rho_\tau$ decomposes
into direct sum of two distinct irreducible representations
$\rho_{+}, \rho_-$ of $\A^{\Z_2}, $ and $[\a_{\rho_+}]=
[\a_{\rho_+}]=[\rho_\tau].$ Similarly for an irreducible
representation $\lambda$ of $\A$, we have that if $[\lambda]=[\bar
\lambda],$ then $\lambda$ decomposes as direct sum of two distinct
irreducible representations $\lambda_{+}, \lambda_-$ of $\A^{\Z_2},
$ and $[\a_{\lambda_+}]= [\a_{\lambda_-}]=[\lambda].$ If
$[\lambda]\neq [\bar \lambda],$ then $\lambda$ and $\bar\lambda$
restrict to the same representation (denoted by $\lambda$) of
$\A^{\Z_2}$, and we have $[\a_\lambda]= [\lambda]+[\bar\lambda].$
Denote by $a$ be the number of irreducible $\tau$-twisted solitons
of $\A,$ , $b$ the number of irreducible representations $\lambda$
of $\A$ such that $\lambda=\bar \lambda$, and $c$ the number of
irreducible representations $\lambda$ of $\A$ such that $\lambda\neq
\bar \lambda.$ By Th. 4.16 of \cite{BEK}, we have
$$
a+b+c = 2 b + c
$$
Hence $a=b.$ Note that by Lemma \ref{tirrep} the number of
irreducible representations of $\L_\tau SU(n)$ at level $k$ is $b.$
By Lemma \ref{tt} it follows that Cor. \ref{irrep} gives all
irreducible $\tau$-twisted  representations of $\A.$ We summarize
these results in the following: \btheor\label{allirrep} (1)
$\rho_\tau$ as given in Cor. \ref{irrep} gives all the irreducible
$\tau$-twisted representations of $\A.$ These representations are in
one to one correspondence with irreducible representations of
$\L_\tau SU(n)$ at level $k$; (2) The list of all irreducible
representations of $\A^{\Z_2}$ are as follows:\par For an
irreducible representation $\lambda$ of $\A$, we have that if
$[\lambda]=[\bar \lambda],$ then $\lambda$ decomposes as direct sum
of two distinct irreducible representations $\lambda_{+}, \lambda_-$
of $\A^{\Z_2}; $ If $[\lambda]\neq [\bar \lambda],$ then $\lambda$
and $\bar\lambda$ restrict to the same representation (denoted by
$\lambda$) of $\A^{\Z_2};$ $\rho_+, \rho_-$ where $\rho_\tau$ corresponds
to irreducible representations of $\L_\tau SU(n)$ at level $k.$
\etheor By Th. \ref{Th.2.6}  and Th. \ref{allirrep}, we have proved
the following: \bcoro\label{mtc} The list of irreducible
representations of $\A^{\Z_2}$ as in (2) of Th. \ref{allirrep} give
rise to a unitary modular tensor category as defined in II.5 of
\cite{Tu} by the construction as given in \S1.7 of \cite{X3m}.
\ecoro
\section{Examples of fusion rules}
Th. \ref{allirrep} and Cor.\ref{mtc} give strong constraints on the
fusion rules related to $\A^{\Z_2}.$ In this section we give some
examples of fusion rules by using the results of \S2. The ideas are
similar to section 9 of \cite{KLX}.\par Denote by $N_{\lambda
\r_1}^{\r_2}= \lan \lambda \rho_1, \rho_2 \ran$ where $\lambda$ is
an irreducible representation of $\A$, and $\rho_1,\rho_2$ are
irreducible $\t$-twisted representations of $\A.$

\blemma\label{generalf} There exists a complex valued matrix
$\psi_{\rho}^{(\mu)}$ where $\rho$ denotes irreducible $\t$-twisted
representations of $\A,$ and $\mu=\bar \mu$ labels irreducible
representation of $\A,$ such that
$$
N_{\lambda \r_1}^{\r_2}= \sum_{\mu, \mu=\bar{\mu}}
\frac{S_{\lambda\mu}}{S_{1\mu}}
\psi_{\r_1}^{(\mu)}\psi_{\r_2}^{(\mu)*}.
$$ \elemma \proof First we assume that $\lambda=\bar\lambda.$ We have
$$
\lan \lambda\rho_1,\rho_2\ran =\lan \lambda \a_{\rho_{1+}},
\a_{\rho_{2+}}\ran =\lan (\lambda_+ +\lambda_-)\rho_{1+},
{\rho_{2+}}\ran
$$

Using Verlinde formula we have:

$$
 \lan (\la_+ +\la_-)  {\rho_{1+}},
{\rho_{2+}}\ran = \sum_{\dot \mu}
(\frac{S_{\la_+\dot\mu}S_{\rho_{1+}\dot\mu}S_{\rho_{2+}\dot\mu}^*}{S_{1\dot\mu}^3}
+\frac{S_{\la_-\dot\mu}S_{\rho_{1+}\dot\mu}S_{\rho_{2+}\dot\mu}^*}{S_{1\dot\mu}^3})
$$
By Lemma 9.1 of \cite{KLX}, we have that $S_{\rho_{+}\dot\mu}=0$ if
$\dot\mu$ comes from restriction of representation $\mu\neq \bar\mu$
of $\A $, and $S_{\la_+ \dot\mu}= - S_{\la_- \dot\mu}$ if $\dot\mu$
comes from restriction of a $\tau$-twisted representation of $\A.$
If $\mu=\bar{\mu},$  by (4) of Lemma 9.1 in \cite{KLX} we have
$$
\frac{S_{\lambda_+\mu_+}S_{\rho_{1+}\mu_+}S_{\rho_{2+}\mu_+}^*}{S_{1\mu_+}^3}
=\frac{S_{\lambda_-\mu_-}S_{\rho_{1+}\mu_-}S_{\rho_{2+}\mu_-}^*}{S_{1\mu_-}^3}
=\frac{S_{\lambda\mu}S_{\rho_{1+}\mu_+}S_{\rho_{2+}\mu_+}^*}{S_{1\mu}S_{1\mu_+}^2}
$$
Set $\psi_{\rho}^{(\mu)}=
\sqrt{2}\frac{S_{\rho_{+}\mu_+}}{S_{1\mu_+}} $ and the Lemma
follows. The case when $\lambda\neq \bar\lambda$ is similar.
\endproof
\subsection{n=3 case}
Lemma \ref{generalf} determines the spectrum of square matrix
$N_\lambda$ whose entries are non-negative integers.  In this
section we determine (up to permutation) $N_\lambda$ for the first
non-trivial case $n=3$. Our results agree with the ansatz given by
\cite{PZ} based on heuristic arguments. \blemma\label{pf}
$d=(d_\rho)_\rho$ is the unique (up to scalar multiplication)
Perron-Frobenius eigenvector of $N_v$ with eigenvalue
$S_{v1}/S_{11}$ where $v$ denotes the vector representation.\elemma
\proof We note that the matrix $N_v$ is
irreducible. In fact since $\rho_1\bar\rho_2\succ \lambda$ for some
$\lambda$, and $\lambda\prec v^m$ for some integer $m$, it follows
that the $(\rho_1,\rho_2)$-th entry of $N_v$ is positive. Hence by
\cite{Gan} the lemma is proved. \endproof\par

We shall refer to the equation
\begin{equation}\label{pfe}
N_v d =\frac{S_{v1}}{S_{11}} d
\end{equation}
as Perron-Frobenius equation.

 Note that by Lemma \ref{generalf}, $N_v= N_{\bar v}.$ Since
every irrep $\la$ can be written as polynomials in $v, \bar v$, it
is sufficient to determine $N_v.$

Let $M=N_v-I.$ By Lemma \ref{generalf} and equation (\ref{kacp}) the
spectrum of $M$ is given as follows.
If $k=2m-1$ is odd, then the spectrum of $M$ is given by $2\cos
(\frac{\pi{(i+1)}}{m+1}), 0\leq i\leq m-1,$ and it is the same as
the fusion matrix $N_{1/2}$ associated with $SU(2)_{m-1}$ where
$1/2$ denotes the spin $1/2$ representation. If $k=2m-2$ is even,
then the spectrum of $N_v$ is given by $2\cos
(\frac{2\pi{(i+1)}}{2m+1}), 0\leq i\leq m-1,$ and it is the same as
the fusion matrix $N_{1}$ associated with $SU(2)_{2m-1}$ acting on
the set of integer spin representations of $SU(2)_{2m-1},$ where $1$
denotes the spin $1$ representation of $SU(2)_{2m-1}.$

Our goal in this section is to show that up to permutation
$N_v=N_{1/2}+I$ when $k=2m-1$ and $N_v=N_{1}$ when $k=2m-2.$\par

First note that since $||M||<2,$ the entries of $M$ can take only
$1,-1,0,$ and since $M=N_v-I,$ only diagonal entries on $M$ can be
$-1.$  By the known spectrum of $M$ we have $tr(M)=0, tr(M^2)=2m-2$
when $k=2m-1$, and $tr(M)=-1, tr(M^2)=2m-1$ when $k=2m-2.$

 Also since $||M^2||<4$ each row of $M$ contains
at most three nonzero entries.

Denote by $k_1,k_2,k_3$ respectively the number of rows of $M$ with
one, two, three nonzero entries respectively. Then we have
$k_1+k_2+k_3=n, k_1+2k_2+3k_3= tr(M^2)= 2m-2$ when $k=2m-1$, and
$k_1+k_2+k_3=n, k_1+2k_2+3k_3= tr(M^2)= 2m-1$ when $k=2m-2.$

Hence $k_1=k_3+2$ when $k$ is odd and $k_1=k_3+1$ when $k$ is even.
For simplicity we enumerate the $\tau$-twisted solitons by $1,...,m$
We associate a graph $G$ to these solitons with vertices $1,...,m$
and the connect $i$ and $j$ ($i\neq j$) by the $(i,j)$-th entry of
$M.$ By Lemma \ref{pf} $G$ is connected.\par

If $k_3=0,$ then $k_1=2$ or $k_1=1$ depending on if $k$ is odd or
even. Permute solitons if necessary, we may assume that the first
row contains only one nonzero entry. Using the equation (\ref{pfe})
we must have $1$ in the first row of $M$, and it is not on the
diagonal. Assume that $2$ is the vertex connected to $1$ on $G$ and
use  the equation \ref{pfe} and the fact that $G$ is connected, we
conclude that unless $m=2,$ $2$ is connected to a new vertex $3.$
Continue this argument we have shown that up to permutation
$N_v=N_{1/2}+I$ when $k$ is odd, and $N_v=N_1$ when $k$ is even.

When $k$ is odd and if all diagonal entries of $M$ are zeros, then
$2I-M$ is a positive definite matrix with all diagonal entries equal
to $2.$ Hence we can find a basis $\{ \epsilon_1, ... \epsilon_m\}$
in $\R^m$ such that $(\epsilon_i,\epsilon_j)=0$ or $-1$ if $i\neq j$
and $(\epsilon_i,\epsilon_i)=2, 1\leq i\leq m.$ It follows that $G$
is a connected Coexter graph , and $G$ must be $A-D-E$ graph (cf.
Page 60 of \cite{Hum} or \S1.4 of \cite{GHJ}). Since $G$ has norm
$2\cos(\frac{\pi}{m+1})$ with $m$ vertices, by inspecting table
1.4.5 of \cite{GHJ} we conclude that up to permutation
$M=N_{1/2}.$\par When $k$ is even, $tr(M)=-1,$ there is at least one
$-1$ on the diagonal of $M$. \par

For the rest of this section we assume that $k_3>0$ and there is at
least one nonzero entry on the diagonal of $M.$ Since one can easily
determine $M$ for $m\leq 3$, we will also assume that $m\geq 4.$ We
will derive contradictions from these assumptions.\par

The basic idea is contained in Remark (2) on Page 23 of \cite{GHJ}.
Introduce new rows numbered by $1',2'...,m'$ and we use $M_1$ to
denote a symmetric $2m\times 2m$ matrix whose $(i,j')$-th entry is
the $(i,j)$-th entry of $M,$ and all other entries of $M_1$ are
equal to zero. We associate a graph $G_1$ to $M_1$ whose vertices
are $1,2,...,m,1',2',...,m'$ and $i$, $j'$ are connected by the
absolute value of the $(i,j)$-th entry of $M.$ Since $G$ is
connected, and by our assumption there is at least one nonzero entry
on the diagonal of $M,$ it follows that $G_1$ is connected. Let
$P=2I-M_1.$ Note that $P$ is positive definite, and we can find a
basis $\e_1,...\e_m, \e_1',...\e_m'$ in $\R^{2m}$ such that the
inner product matrix of this basis is $P.$ By definition
$\e_1,...,\e_m$ and $\e_1',...,\e_m'$ are two orthogonal sets, and
the angles between the lines spanned by the elements in the basis
are either $\pi/3$ or $\pi/2,$ and by Th. 3.5 of \cite{CGS} we
conclude that $\e_1,...\e_m, \e_1',...\e_m'$ is contained in a
direct sum of root systems of $A-D-E.$ Since $G_1$ is connected, we
conclude that $\e_1,...\e_m, \e_1',...\e_m'$ is contained in one
root system. If this root system is $E_6$ or $E_8$, then we have
$m=3$ or $m=4$. One can easily rule out these two cases using
$k_1=k_3+2$ or $k_1=k_3+1$ and equation (\ref{pfe}).  If this root
system is $A_{2m+1},$ since the elements of this system are of the
form $e_i-e_j, 1\leq i,j\leq 2m+1$ (cf. Definition 3.1 of
\cite{CGS}), it follows that $M$ can not have a row with three
nonzero entries. Hence to finish the proof we assume that
$\e_1,...\e_m, \e_1',...\e_m'$ is contained in  root system of
$D_{2m}.$
Root system of type $D_{2m}$ consists of vectors $\pm e_i\pm e_j,
1\leq i\neq j\leq 2m,$ where $e_i, 1\leq i\leq 2m$ is an orthonormal
basis in $\R^{2m}.$ By assumption we can identify $\e_1,...\e_m,
\e_1',...\e_m'$ as a subset of $\pm e_i\pm e_j, 1\leq i\neq j\leq
2m.$

The following lemma will be used repeatedly in the following, and
its proof follows directly from definitions: \blemma\label{pair} (1)
If a vertex on $G_1$ is connected to three different vertices
$i_1',i_2',i_3'$ then $\e_{i_1}',\e_{i_2}',\e_{i_3}'$ must contain
(up to multiplication by $-1$) $e_j+e_k, e_j-e_k$ for some $1\leq
j\neq k\leq 2m;$\par (2) If two vertices $i',j'$ of $G_1$ are such
that $\e_i'=e_j+e_k, \e_j'=e_j-e_k,$ then $i'$ is connected to
vertex $l$ in $G_1$ iff $j'$ is connected to vertex $l;$ \par

(3) If $\e_{i_1}'=e_j+e_k, \e_{i_2}'=e_j-e_k, \e_{j_1}=e_p+e_q,
\e_{j_2}=e_p-e_q,$ then $\{j,k \} \cap \{p, q\} = \emptyset.$

\elemma
\bdefin\label{type} A vertex $i$ of $G$ is called type $0$ if the
$(i,i)$-the entry of $M$ is zero, and type $1$ otherwise. \edefin

Case (1): If a type $1$ vertex is connected two different type $0$
vertices on $G$ , by Lemma \ref{pair} the two type $0$ vertices are
not connected.  By permuting basis elements if necessary, we may
assume that $(\e_i,\e_j')=0, |(\e_i,\e_i')|=1, 1\leq i\neq j\leq 3.$
Since $\e_1,...\e_m, \e_1',...\e_m'$ is contained in  root system of
$D_{2m},$ $\e_1',...\e_m'$ and $\e_1,...,\e_m$ are two orthogonal
sets, it follows that we may assume that (up to multiplication by
$-1$) that $\e_1'=e_3+e_5, \e_2'=e_1+e_2, \e_3'=e_1-e_2$ and
$\e_1=e_1+e_5, \e_2=e_3+e_4, \e_3=e_3-e_4.$ But then $e_5$ is in the
subspace spanned by $\e_1',\e_2,\e_3$ and also in the subspace
spanned by $\e_1, \e_2',\e_3'$, contradicting the fact that
$\e_1,...\e_m, \e_1',...\e_m'$ is a basis.\par Case (2): If two type
$1$ vertices are connected on $G$, assume that one such vertex is
connected to either a type $0$ or type $1$ vertex. The first case is
ruled out by Lemma \ref{pair} and Case (1), and the second case is
impossible by Lemma \ref{pair} and the fact that  $G$ is connected
and $m\geq 4.$\par Case (3): If two type $1$ vertices are connected
to the same type $0$ vertex on $G$, since $m\geq 4$ and $G$ is
connected, by case (1) and (2) we assume that the type $0$ vertex is
connected to one additional vertex on $G,$ but one checks easily
that this is impossible by Lemma \ref{pair}. \par Now consider a
subgraph $G'$ obtained from $G$ by deleting type $1$ vertices and
edges with one endpoint a type $1$ vertex. Since $G$ is connected,
by Case (1) $G'$ is connected. Moreover, if $i$ is a type $0$
vertex, define $\eta_i: =\frac{1}{\sqrt 2}(\e_i+\e_i')$. Then
$(\eta_i,\eta_i)=2$,  and it follows that two vertices $i_1,i_2$ on
$G'$ are connected by $-(\eta_i,\eta_j)$ edges. Hence $G'$ is a
connected Coexter graph, and by Page 60 of \cite{Hum} we know that
$G'$ is an $A-D-E$ graph. Since we assume that $k_3\geq 1$, $G'$
must be type $D$ or $E.$ In the case $k_1=k_3+2$, we must have no
type $1$ vertex attached to the end points of $G'$, contracting our
assumption that type $1$ vertex exists and must be connected to one
type $0$ vertex by case (2).  In the case $k_1=k_3+1$, we must have
exactly one type $1$ vertex connected to one endpoint of either type
$D$ or type $E$ graph, and these cases can be directly  ruled out by
tedious calculations using the Perron-Frobenius equation. Here we
give a different approach. By fusion rules of $SU(2)_{2m-1}$ in
equation (\ref{su2}) , $N_{m-1}$ can be written as a polynomial of
$N_1$ with integer coefficients, and we have $N_{m-1}^2= N_1 +I.$
The spectrum of $N_{m-1}$ is given by $2\cos(\frac{\pi(i+1)}{2m+1}),
0\leq i\leq m-1.$ Since $N_1$ has the same spectrum as $N_v$, it
follows that there is symmetric matrix $M'$ with integer entries
such that $M'^2= N_v+I,$ and $M'$ has the same spectrum as
$N_{m-1}.$ Since $N_v$ has one $1$ and $m-1$ $2$'s on the diagonal,
it follows that $M'$ has one row with one nonzero entry $\pm 1$, and
$m-1$ rows with two nonzero entries which are $\pm 1.$ Now associate
a graph with vertices $1,...,m$ to $M'$ so that the $i$-th and the
$j$-th vertex are connected by the absolute value of the $(i,j)$-th
entry of $M'$ (If $i=j$ and the $(i,i)$-th entry of $M'$ is $\pm 1$,
we connect $i$ to itself by a loop). As in the case with $M$ this
graph is connected, and it follows that it is a line segment with
one loop attached to an endpoint. Since $G'$ has a trivalent vertex,
it follows that there exists four different vertices $i,j,k,l$ such
that the $(i,j), (i,k), (i,l)$ entries of $N_v$ are $1,$ but one
checks easily that this is impossible since $M'^2= N_v+I.$ We have
proved the following:
\bprop\label{n3fusion} (1) When $k=2m-1$, we can label
$\tau$-twisted irrep of $\A$ by integers $1,2,...,m$ such that the
$(i,j)$-th entry of $N_v$ is given by $\delta_{ij}+ \lan \frac{1}{2}
\frac{(i-1)}{2},\frac{(j-1)}{2}\ran$ where half integers $0\leq
k/2\leq {m-1}/2$ label the spin of irreps of $SU(2)_{m-1};$\par (2)
When $k=2m-2$, we can label  $\tau$-twisted irrep of $\A$ by
integers $1,2,...,m$ such that the $(i,j)$-th entry of $N_v$ is
given by $\delta_{ij}+ \lan 1 (i-1),(j-1)\ran$ where  integers
$0\leq k\leq {m-1}$ label the spin of irreps of $SU(2)_{2m-1}.$
\par\eprop
By Prop. \ref{n3fusion} the Perron-Frobenius eigenvector $d$ is up
to multiplication by a  positive constant $\delta$ and possible
permutations equal to $(d_{(i-1)/2})_{1\leq i\leq m}$ (resp.
$(d_{i-1})_{1\leq i\leq m})$) when $k=2m-1$ (resp. when $k=2m-2$).
By Prop. 3.1 of \cite{BEK} $\sum_{\rho_\tau} d_{\rho_\tau}^2 =
\sum_\lambda d_\lambda^2,$ and by equation (\ref{kacp}) we can
determine $\delta$ uniquely. When $k=2m-1$, $\delta=
\frac{\sqrt{3(m+1)}}{2\sin^2(\frac{\pi}{2m+2})},$ and when $k=2m-2,$
 $\delta=
\frac{\sqrt{3(2m+1)}}{4\sin(\frac{\pi}{2m+1})\sin(\frac{2\pi}{2m+1})}.$
We have therefore proved the following \bcoro\label{n3index} (1)
When $k=2m-1$,  the set of indices of $\tau$-twisted solitons are
given by $\{ \frac{3(m+1) }{4\sin^4(\frac{\pi}{2m+2})}
\frac{\sin(i\pi/(m+1))}{\sin(\pi/(m+1))}, 1\leq i\leq m \};$ \par
(2)
 When
$k=2m-2$,  the set of indices of $\tau$-twisted solitons are given
by $\{
\frac{3(2m+1)}{16\sin^2(\frac{\pi}{2m+1})\sin^2(\frac{2\pi}{2m+1}) }
\frac{\sin((2i-1)\pi/(2m+1))}{\sin(\pi/(2m+1))}, 1\leq i\leq m \}.$
\ecoro

\end{document}